\documentclass[smallextended,envcountsect]{svjour3}
\usepackage{graphicx,color,cite}
\usepackage{mathrsfs}
\usepackage{amsmath, amscd, amsfonts, amssymb, graphicx, color}
\usepackage[bookmarksnumbered, colorlinks, plainpages]{hyperref}
\usepackage{mathptmx}   
\usepackage[colorinlistoftodos]{todonotes}
\usepackage[left=4cm, right=4cm, top=3.5cm, bottom=3.5cm]{geometry}

\usepackage{latexsym}

\if{
\newtheorem{theorem}{Theorem}[section]

\newtheorem{lemma}{Lemma}[section]

\newtheorem{remark}{Remark}[section]

\newtheorem{claim}{\textbf{Claim}}[section]

\newtheorem{corollary}{Corollary}[section]
\newtheorem{proposition}{Proposition}[section]

}\fi

\renewcommand\theenumi{\roman{enumi}}

\newcounter{mycount}

\smartqed

\usepackage{etex}

\usepackage{mathtools}

\usepackage{enumitem}
\renewcommand\theenumi{(\roman{enumi})}
\renewcommand\theenumii{(\alph{enumii})}

\usepackage{csquotes}

\usepackage{todonotes}

\makeatletter
\let\orgdescriptionlabel\descriptionlabel
\renewcommand*{\descriptionlabel}[1]{
 \let\orglabel\label
 \let\label\@gobble
 \phantomsection
 \edef\@currentlabel{#1}
 \let\label\orglabel
 \orgdescriptionlabel{#1}
}
\def\th@plain{
 \thm@notefont{}
 \itshape
}
\def\th@definition{
 \thm@notefont{}
 \normalfont
}

\g@addto@macro\th@definition{\thm@headpunct{}}
\g@addto@macro\th@plain{\thm@headpunct{}}
\makeatother

\usepackage{subfig}

\usepackage[final]{showkeys}

\usepackage{etoolbox}

\usepackage{mathrsfs}

\usepackage[titletoc]{appendix}
\usepackage[doc]{optional}

\usepackage{soul}

\usepackage{colortbl,booktabs,
multirow}
\usepackage{xcolor}

\usepackage{cancel}

\usepackage{empheq}

\definecolor{myblue}{rgb}{.8, .8, 1}

\usepackage{hyperref}
\hypersetup{
colorlinks=true,
linkcolor=blue,
citecolor=red,
filecolor=magenta,
urlcolor=cyan
}

\usepackage[
open,
openlevel=2,
atend,
numbered
]{bookmark}

\usepackage[capitalize, nameinlink, noabbrev]{cleveref}
\crefname{equation}{}{}
\crefname{chapter}{Chapter}{Chapters}
\crefname{item}{item}{items}
\crefname{figure}{Figure}{Figures}
\crefname{theorem}{Theorem}{Theorems}
\crefname{lemma}{Lemma}{Lemmas}
\crefname{proposition}{Proposition}{Propositions}
\crefname{corollary}{Corollary}{Corollarys}
\crefname{definition}{Definition}{Definitions}
\crefname{fact}{Fact}{Facts}
\crefname{example}{Example}{Examples}
\crefname{algorithm}{Algorithm}{Algorithms}
\crefname{remark}{Remark}{Remarks}
\crefname{note}{Note}{Notes}
\crefname{notation}{Notation}{Notations}
\crefname{case}{Case}{Cases}
\crefname{exercise}{Exercise}{Exercises}
\crefname{question}{Question}{Questions}
\crefname{claim}{Claim}{Claims}
\crefname{enumi}{}{}

\usepackage{pgf}

\usepackage{array}
\usepackage{tabu}

\allowdisplaybreaks

\numberwithin{equation}{section}

\renewcommand\theenumi{(\roman{enumi})}
\renewcommand\theenumii{(\alph{enumii})}

\spnewtheorem*{Proof}{Proof.}{\bf}{\rm}

\begin{document}

\title{{Complexity of Error Bounds for Systems of Linear Inequalities}  \thanks{This work was supported by the Natural Science Foundation of Hebei Province (A2024201015) and the Excellent Youth Research Innovation Team of Hebei University (QNTD202414).
}}

\titlerunning{{Complexity of Error Bounds for Systems of Linear Inequalities}}

\author{Zhou Wei  \and Michel  Th\'era \and Jen-Chih Yao}

\institute{Zhou Wei\at Hebei Key Laboratory of Machine Learning and Computational Intelligence \& College of Mathematics and Information Science, Hebei University, Baoding, 071002, China\\ 	ORCID 0000-0001-9697-7198\\
	\email{weizhou@hbu.edu.cn}\\
	Michel Th\'era \at XLIM UMR-CNRS 7252, Universit\'e de Limoges, Limoges, France\ \\ORCID 0000-0001-9022-6406 \\ \email{michel.thera@unilim.fr}\\
	Jen-Chih Yao \at Research Center for Interneural Computing, China Medical University Hospital,
	China Medical University, Taichung, Taiwan \\ and 
	\at Academy of Romanian Scientists, 50044 Bucharest, Romania\\ \email{yaojc@mail.cmu.edu.tw}
}

\date{Received: date / Accepted: date}

\maketitle

\begin{abstract}

Error bounds have been studied for more than seventy years, beginning with the seminal result of Hoffman (1952) [{\it J. Res. Natl. Bur. Standards}, 49 (1952), 263--265], which establishes an upper bound for the distance from an arbitrary point to the solution set of a linear system. Despite this long history, relatively little is known about the intrinsic computational complexity of error bounds.

In this paper, we investigate the complexity of error bounds for systems of linear inequalities. We first reformulate the problem as a finite collection of min--max optimization problems defined on the unit sphere and associated with subsets of the rows of the given matrix. We then prove that the problem does not belong to the class {\bf P}, while it is {\bf co\mbox{-}NP}-complete.

Furthermore, we establish the existence of a pseudo-polynomial-time algorithm for solving the complementary problem. In particular, the complement may be regarded as a number problem, although it is not {\bf NP}-complete in the strong sense unless {\bf P} = {\bf NP}.


\end{abstract}

\section{Introduction}

It is well known that the study of error bounds can be traced back to the
pioneering work of Hoffman in 1952. More precisely, given an
$m\times n$ matrix $A$ and an $m$-vector $b$, Hoffman
\cite{Hoffman1952} established a fundamental estimate for the Euclidean
distance from a point $x$ to the polyhedral set
\[
P_{A,b}:=\{u\in\mathbb{R}^n: Au\leq b\}.
\]
Specifically, he proved the existence of a constant $c>0$ such that
\begin{equation}\label{1.1}
{\bf d}(x,P_{A,b})\leq c\|(Ax-b)_+\|
\end{equation}
for every $x\in\mathbb{R}^n$, where $(Ax-b)_+$ denotes the positive part of
the vector $(Ax-b)$, namely,
\[
(z)_+:=\max\{z,0\}
\]
for each component $z$ of $(Ax-b)$.

Inequality \eqref{1.1} shows that the Euclidean distance from a point
$x$ to the feasible set $P_{A,b}$ can be estimated above by a constant
multiple of the residual error $\|(Ax-b)_+\|$. Hoffman's theorem has become a cornerstone of variational analysis and optimization and
has been extensively studied and generalized by numerous authors,
including Robinson~\cite{Robinson1975},
Mangasarian~\cite{Mangasarian1985},
Auslender and Crouzeix~\cite{AC1988},
Khachiyan~\cite{Khachiyan1995},
Tuncel~\cite{Tuncel1999},
Pang~\cite{Pang1997},
Lewis and Pang~\cite{LewisPang1996},
Klatte and Li~\cite{KlatteLi1999},
Jourani~\cite{jourani},
Abassi and Th\'era~\cite{Abassi-Thera1,Abassi-Thera},
Ioffe~\cite{ioffe-book},
Mordukhovich~\cite{Boris1,Boris2},
Penot~\cite{penot-book},
Thibault~\cite{thibault-book},
among many others.

The theory of error bounds has remained an active area of research,
particularly in the development of necessary and sufficient conditions
for their existence, as well as efficient methods for their computation.
Error bounds are now well established in optimization and variational
analysis and have numerous important applications, including sensitivity
analysis for linear and integer programming problems
(see \cite{Rob73,Rob77}) and convergence analysis of descent methods for
linearly constrained minimization problems
(see \cite{HLu,TsB93}). They also play a crucial role in feasibility
problems involving intersections of finitely many convex or nonconvex
closed sets (see \cite{5,6,7,BurDeng02}) and have found important
applications in image reconstruction (see \cite{16}).

Furthermore, error bounds are closely connected with several central
concepts in variational analysis, including weak sharp minima
(cf. \cite{BD2}), metric regularity and metric subregularity of
set-valued mappings
(cf. \cite{ioffe-JAMS-1,ioffe-JAMS-2,Kru15,Kru15.2}),
the Aubin property of set-valued mappings
(cf. \cite{CKLT,Kruger-LY}),
and transversality and subtransversality of collections of closed sets
(see \cite{KLT2017,WZ,WTY2024-SVVA}).

The study of error bounds relies on a variety of mathematical tools and
concepts, including derivatives, subdifferentials, normal cones, and
other generalized differential constructions. Broadly speaking, two
main approaches have been developed for deriving sufficient and/or
necessary conditions for error bounds: one based on points belonging to
the solution set and the other based on points outside the solution set.
For further details on these approaches and related criteria, we refer
the reader to Ioffe \cite{Ioffe1979},
Lewis and Pang \cite{LewisPang1996},
Pang \cite{Pang1997},
Z\u{a}linescu \cite{Za2001},
Wu and Ye \cite{55},
Zheng and Ng \cite{ZN2004},
Ngai and Th\'era \cite{NT2004,Huynh2008Error,NT2009},
Fabian, Henrion, Kruger, and Outrata \cite{FHKO2010},
Li, Meng, and Yang \cite{LMY2018},
Dutta and Mart\'inez-Legaz \cite{DuMa2021},
Cuong and Kruger \cite{CK2022},
and Wei, Th\'era, and Yao \cite{WTY2023}.

Given the fundamental role played by error bounds in variational
analysis and optimization, one may naturally expect such properties to
be readily verifiable in practical applications. This observation leads
to the following natural question: how difficult is it, from a
computational viewpoint, to determine whether a system admits an error
bound? Motivated by this question, we investigate in this paper the
computational complexity of error bounds for systems of linear
inequalities. More specifically, we examine whether such bounds can be
computed in polynomial time and analyze the complexity of the associated
decision problems.

Our aim is to study the complexity of error bounds for systems of linear
inequalities. Given an $m\times n$ matrix $A$ and an $m$-vector $b$, we
show that the problem of determining error bounds for the system
$Ax\leq b$ can be reformulated equivalently as solving finitely many
min--max optimization problems over the unit sphere associated with
certain subsets of rows of $A$ (see \cref{pro3.1}). We then prove that
there is no polynomial-time algorithm for solving the error bound
problem, since the number of min--max problems that must be verified to
have negative optimal values may grow exponentially in the worst case
(see \cref{th3.1}). Next, we establish that the problem is
{\bf co\mbox{-}NP}-complete and hence {\bf NP}-hard
(see \cref{th3.2} and \cref{coro3.3}). Furthermore, we prove the
existence of a pseudo-polynomial-time algorithm for solving the
complementary problem (see \cref{th3.3}). In particular, the complement
can be viewed as a number problem and is not {\bf NP}-complete in
the strong sense unless ${\bf P}={\bf NP}$.

The remainder of the paper is organized as follows. In Section~2, we
present the definitions and preliminary results used throughout the
paper. Section~3 is devoted to the complexity analysis of error bounds
for systems of linear inequalities of the form $Ax\leq b$, where $A$ is
a given $m\times n$ matrix and $b$ is an $m$-vector. Finally,
concluding remarks and perspectives for future research are provided in
 Section~4.

\section{Preliminaries}

We recall some preliminaries used in the sequel, beginning with the notion of a problem.

\begin{itemize}
\item
A {\it problem} is a general question depending on parameters whose values are initially unspecified. It is specified by (i) the definition of its parameters and (ii) the properties required of any feasible solution. An {\it instance} is obtained by fixing all parameters.

\item
A function \( f : S \to \mathbb{R}_+ \) is {\it polynomially bounded} by \( g : S \to \mathbb{R}_+ \) if there exists a polynomial \( \phi \) such that
\[
f(s) \le \phi(g(s)) \quad \forall s \in S.
\]

\item
An {\it algorithm} solves a problem if it produces a correct solution for every instance in finite time. It runs in {\it polynomial time} if its running time (number of arithmetic operations, denoted \( f \)) is polynomially bounded by the input size \( g \). Such an algorithm is called a {\it polynomial-time algorithm}.

\item
A problem is in {\bf P} if it admits a polynomial-time algorithm. Linear programming is in {\bf P} due to Karmarkar's algorithm.

\item
A {\it decision problem} is a problem whose answer is either "yes" or "no". The class {\bf NP} consists of decision problems for which "yes" instances admit polynomially verifiable certificates. The class Co-{\bf NP} is defined analogously for "no" instances. For example, mixed-integer linear programming feasibility is in {\bf NP} (cf. \cite[Theorem 4.36]{CCZ2014}).

\item
A problem \( Q \in {\bf NP} \) is {\it {\bf NP}-complete} if every \( D \in {\bf NP} \) reduces to \( Q \) in polynomial time. SAT is {\bf NP}-complete \cite{Cook71}, and so is feasibility of mixed-integer linear sets.

\item
A problem \( Q \in \text{\bf Co-NP} \) is {\it Co-{\bf NP}-complete} if every \( D \in \text{\bf Co-NP} \) reduces to \( Q \) in polynomial time.

\item
A (not necessarily decision) problem \( Q \) is {\bf NP}-{\it hard} if every \( D \in {\bf NP} \) reduces to \( Q \) in polynomial time.
\end{itemize}

\medskip

Let \( \mathbb{R}^n \) be equipped with inner product \( \langle x,y\rangle := x^T y \) and norm \( \|x\| := \sqrt{\langle x,x\rangle} \). Let \( {\bf B}_{\mathbb{R}^n} \) denote the closed unit ball.

For finitely many points \( d_1,\dots,d_k \in \mathbb{R}^n \), let \( \mathrm{conv}(d_1,\dots,d_k) \) denote their convex hull. For \( D \subset \mathbb{R}^n \), let \( \mathrm{int}(D) \) and \( \mathrm{bdry}(D) \) denote its interior and boundary, and define
\[
\mathbf{d}(x,D) := \inf_{u \in D} \|x-u\|,
\quad \mathbf{d}(x,\emptyset):=+\infty.
\]

Let \( \Gamma_0(\mathbb{R}^n) \) denote the class of proper, lower semicontinuous convex functions \( \varphi : \mathbb{R}^n \to \mathbb{R}\cup\{+\infty\} \), with \( \mathrm{dom}(\varphi) := \{x : \varphi(x) < +\infty\} \neq \emptyset \).

For \( \varphi \in \Gamma_0(\mathbb{R}^n) \), \( \bar{x} \in \mathrm{dom}(\varphi) \), and \( h \in \mathbb{R}^n \), the directional derivative is
\[
d^+\varphi(\bar{x},h)
:= \lim_{t \downarrow 0} \frac{\varphi(\bar{x}+th)-\varphi(\bar{x})}{t}.
\]
The function \( t \mapsto \frac{\varphi(\bar{x}+th)-\varphi(\bar{x})}{t} \) is nonincreasing as \( t \downarrow 0 \) \cite{Ph}, and thus
\[
d^+\varphi(\bar{x},h)
= \inf_{t>0} \frac{\varphi(\bar{x}+th)-\varphi(\bar{x})}{t}.
\]

The subdifferential of \( \varphi \) at \( \bar{x} \) is
\[
\partial \varphi(\bar{x})
:= \{\xi \in \mathbb{R}^n : \xi^T(x-\bar{x}) \le \varphi(x)-\varphi(\bar{x}),\ \forall x \in \mathbb{R}^n\}.
\]

If \( \varphi \) is continuous at \( \bar{x} \), then \( \partial \varphi(\bar{x}) \neq \emptyset \), and
\[
\partial \varphi(\bar{x})
= \{\xi \in \mathbb{R}^n : \xi^T h \le d^+\varphi(\bar{x},h),\ \forall h \in \mathbb{R}^n\},
\]
\[
d^+\varphi(\bar{x},h)
= \max_{\xi \in \partial \varphi(\bar{x})} \xi^T h.
\]

\setcounter{equation}{0}

\section{Main Results}

In this section, we investigate the computational complexity of error bounds for systems of linear inequalities. We begin with several auxiliary results that play a fundamental role in the subsequent complexity analysis and are also of independent interest.

\begin{lemma}\label{lem3.2}
The problem: " Given a positive integer \( n \) and a finite set \( \{d_1,\ldots,d_k\}\subseteq\mathbb{R}^n \), determine whether
$
\min_{\|h\|=1}\max_{1\le i\le k} d_i^T h<0 "
$
belongs to the class {\bf P}.
\end{lemma}

\begin{proof}
Let \( d_1,\ldots,d_k\in\mathbb{R}^n \), and define
\[
D:=(d_1^T,\ldots,d_k^T), \qquad \mathbf e:=(1,\ldots,1)\in\mathbb{R}^k.
\]
We claim that
\begin{equation}\label{3.1-260223}
\min_{\|h\|=1}\max_{1\le i\le k} d_i^T h<0
\end{equation}
holds if and only if the linear system
\begin{equation}\label{3.9}
\left\{
\begin{aligned}
&\begin{pmatrix}D\\ \mathbf e\end{pmatrix}\lambda
 = \begin{pmatrix}0\\1\end{pmatrix},\\
&0\le \lambda_i\le 1,\quad i=1,\ldots,k,\\
&\lambda=(\lambda_1,\ldots,\lambda_k)\in\mathbb{R}^k
\end{aligned}
\right.
\end{equation}
is inconsistent.

Assuming the equivalence, the conclusion follows immediately from the existence of polynomial-time algorithms for linear programming feasibility problems.

Suppose first that \eqref{3.1-260223} holds and let \( \bar\lambda \) be feasible for \eqref{3.9}. Then
\[
\sum_{i=1}^k \bar\lambda_i=1,
\qquad
\sum_{i=1}^k \bar\lambda_i d_i=0.
\]
By \eqref{3.1-260223}, there exists \( \bar h\in\mathbb{R}^n \) with \( \|\bar h\|=1 \) such that
\[
\max_{1\le i\le k} d_i^T\bar h<0.
\]
Hence,
\[
0
=
\Big\langle \sum_{i=1}^k \bar\lambda_i d_i,\bar h\Big\rangle
=
\sum_{i=1}^k \bar\lambda_i d_i^T\bar h
<0,
\]
a contradiction. Therefore, \eqref{3.9} is inconsistent.

Conversely, assume that \eqref{3.9} is inconsistent. Then
\[
0\notin \mathrm{conv}(d_1,\ldots,d_k).
\]
By the strict separation theorem, there exists \( \hat h\in\mathbb{R}^n \) with \( \|\hat h\|=1 \) such that
\[
\max_{1\le i\le k} d_i^T\hat h<0.
\]
Consequently,
\[
\min_{\|h\|=1}\max_{1\le i\le k} d_i^T h
\le
\max_{1\le i\le k} d_i^T\hat h
<0,
\]
which proves \eqref{3.1-260223}.\qed
\end{proof}

\begin{corollary}\label{coro3.1}
The problem 
``Given a positive integer \( n \) and a finite set \( \{d_1,\ldots,d_k\}\subseteq\mathbb{R}^n \), determine whether
$
\min_{\|h\|=1}\max_{1\le i\le k} d_i^T h\ge 0
$ '' 
belongs to the class {\bf P}.
\end{corollary}

The next proposition concerns the complexity of deciding whether
\[
\min_{\|h\|=1}\max_{1\le i\le k} d_i^T h\neq 0
\]
for a finite set \( \{d_1,\ldots,d_k\}\subseteq\mathbb{R}^n \) when the dimension is fixed.

\begin{proposition}\label{lem4.1}
Let \( n \) be fixed. Then the problem
``Given a finite set \( \{d_1,\ldots,d_k\}\subseteq\mathbb{R}^n \), determine whether
$
\min_{\|h\|=1}\max_{1\le i\le k} d_i^T h\neq 0
$ ''
belongs to the class {\bf P}.
\end{proposition}

\begin{proof}
By Lemma~\ref{lem3.2}, it suffices to verify in polynomial time whether
\[
\min_{\|h\|=1}\max_{1\le i\le k} d_i^T h>0.
\]

Let
\[
\Omega:=\mathrm{conv}(d_1,\ldots,d_k).
\]
We first show that
\begin{equation}\label{4.9}
\min_{\|h\|=1}\max_{1\le i\le k} d_i^T h>0
\quad\Longleftrightarrow\quad
0\in \operatorname{int}(\Omega).
\end{equation}

Define
\[
g(x):=\max_{1\le i\le k} d_i^T x,
\qquad x\in\mathbb{R}^n.
\]
Then
\[
\partial g(0)=\Omega,
\qquad
d^+g(0,h)=g(h)
\quad \forall h\in\mathbb{R}^n.
\]
Moreover, by Lemma~\ref{lem3.2} and its proof,
\begin{equation}\label{4.10}
\min_{\|h\|=1}\max_{1\le i\le k} d_i^T h\ge 0
\quad\Longleftrightarrow\quad
0\in\Omega.
\end{equation}
On the other hand, \cite[Remark 6(b)]{WTY2022} yields
\[
\min_{\|h\|=1}\max_{1\le i\le k} d_i^T h=0
\quad\Longleftrightarrow\quad
0\in \operatorname{bdry}(\Omega).
\]
Since
\[
\Omega=\operatorname{int}(\Omega)\cup\operatorname{bdry}(\Omega),
\]
the equivalence in \eqref{4.9} follows from \eqref{4.10}.

Now, by \cite[Theorem 3.5]{Chazelle1993}, the convex hull of \( k \) points in \( \mathbb{R}^n \) can be computed in time
\[
O(k\log k+k^{\lfloor n/2\rfloor}).
\]
Hence, when \( n \) is fixed, \( \Omega \) can be computed in polynomial time and represented as
\[
\Omega=\{x\in\mathbb{R}^n: Mx\le \gamma\},
\]
where \( M \) and \( \gamma \) are computable in polynomial time.

By \cite[Proposition 4.1]{WTY2024},
\[
\operatorname{int}(\Omega)
=
\{x\in\mathbb{R}^n: Mx<\gamma\}.
\]
Therefore,
\[
0\in\operatorname{int}(\Omega)
\quad\Longleftrightarrow\quad
\gamma>0,
\]
which can clearly be verified in polynomial time.

Combining this with \eqref{4.9}, we conclude that deciding whether
\[
\min_{\|h\|=1}\max_{1\le i\le k} d_i^T h>0
\]
is polynomially solvable. The assertion follows.\qed
\end{proof}

\begin{remark}\label{rem4.1}
It follows from \cref{lem4.1} that the problem
``Given a finite set \( \{d_1,\ldots,d_k\}\subseteq\mathbb{R}^n \), determine whether
$
\min_{\|h\|=1}\max_{1\le i\le k} d_i^T h \neq 0
$ ''
is polynomially solvable when the dimension \( n \) is fixed.
It is natural to ask whether this problem remains in {\bf P} when the dimension is not fixed. At present, this question is open. A natural strategy would follow the approach used in the proof of \cref{lem4.1}, which reduces the problem to the computation of the convex hull of a finite point set and to the complexity analysis of convex hull algorithms.

More precisely, the issue reduces to computing a minimal inequality representation
\[
Q=\{x\in\mathbb{R}^n : Mx\le \gamma\}
\]
of a full-dimensional convex polytope \( Q=\operatorname{conv}(S) \), where \( S \) is the input set of points. The relevant complexity measure is the worst-case running time of computing this representation from \( k \) points in \( \mathbb{R}^n \).

For fixed dimension, Chazelle \cite{Chazelle1993} established an optimal convex hull algorithm with worst-case complexity \( O(k^{\lfloor n/2\rfloor}) \) for \( n\ge 4 \). This bound is tight, as shown by the Upper Bound Theorem (see \cite[Theorem 2]{Fukuda2004}), which implies that the worst-case output complexity has the same order of growth.

Despite its theoretical optimality, Chazelle's algorithm is not considered practically efficient in higher dimensions. Under nondegeneracy assumptions, polynomial-time algorithms for convex hull computation are available (cf. \cite{CCH53,Dye83,CK70,AF96}). However, for the general (possibly degenerate) case and variable dimension, no polynomial-time algorithm is known.
\end{remark}

\medskip

{\it Throughout this section, given two positive integers $m,n$, an \( m \times n \) matrix $A$ and an \( m \)-vector  \( b \), we always denote \( A := (a_1^T, \cdots, a_m^T)^T \) and \( b:=(b_1, \cdots, b_m)^T \). We denote the polyhedral set by
	\begin{equation}\label{3.4}
		P_{A,b} := \left\{x \in \mathbb{R}^n : a_i^T x \leq b_i,\ i = 1, \cdots, m\right\}.
\end{equation}
Let $\varphi_{A,b}(x) := \max_{1 \leq i \leq m} (a_i^T x - b_i)$ for each $ x \in \mathbb{R}^n$.
We denote
\begin{equation}\label{3.7}
	J_{A,b}(x) := \left\{ i \in \{1, \cdots, m\} : a_i^T x - b_i = \varphi_{A,b}(x) \right\},
\end{equation}
\begin{equation}\label{3.7-260223}
\mathcal{J}_{A,b} := \left\{ J_{A,b}(x) : \varphi_{A,b}(x) > 0 \right\}.
\end{equation}
}

We aim to investigate the computational complexity of error bounds for systems of linear inequalities. As a main result, we show that the corresponding recognition problem admits no polynomial-time algorithm in the worst case and is, in fact, Co-\textbf{NP}-complete. To this end, we first study global error bounds for general convex functions and establish an equivalent characterization in terms of directional derivatives.

Let \( f \in \Gamma_0(\mathbb{R}^n) \), and define
\begin{equation}\label{3.1}
\mathbf{S}_f := \{x \in \mathbb{R}^n : f(x) \leq 0\},
\end{equation}
the lower level set of \( f \). We say that \( f \) admits a global error bound if there exists \( c \in (0,+\infty) \) such that
\[
c\,\mathbf{d}(x,\mathbf{S}_f) \leq f_+(x), \qquad \forall x \in \mathbb{R}^n,
\]
where \( f_+(x) := \max\{f(x),0\} \). We define the global error bound modulus by
\begin{equation}\label{3.2}
\mathbf{Er}(f) := \inf_{f(x)>0} \frac{f(x)}{\mathbf{d}(x,\mathbf{S}_f)}.
\end{equation}
Thus, \( f \) admits a global error bound if and only if \( \mathbf{Er}(f) > 0 \).

\medskip

The following result, established in \cite[Theorem 3.1]{WTY2025}, provides an exact representation of the global error bound modulus for proper lower semicontinuous convex functions. The proof relies solely on the Ekeland variational principle and plays a central role in our analysis.

\begin{lemma}\label{lem3.1}
Let \( f \in \Gamma_0(\mathbb{R}^n) \) and assume that \( \mathbf{S}_f \neq \emptyset \). Then
\begin{equation}\label{3.3}
\mathbf{Er}(f)
=
\inf_{f(x)>0}\left(-\inf_{\|h\|=1} d^+ f(x,h)\right).
\end{equation}
\end{lemma}

Given an \( m \times n \) matrix \( A \) and a vector \( b \in \mathbb{R}^m \), recall from Hoffman's seminal work \cite{Hoffman1952} that, provided \( P_{A,b} := \{x : Ax \le b\} \neq \emptyset \), the quantity
\begin{equation}\label{3.5}
\inf_{x \notin P_{A,b}}
\frac{\max_{1\le i\le m}(a_i^T x - b_i)_+}{\mathbf{d}(x,P_{A,b})}
\end{equation}
is known as the Hoffman constant of the system \( Ax \le b \), where \( (t)_+ := \max\{t,0\} \).

\medskip

It is straightforward to verify that the Hoffman constant of \( Ax \le b \) coincides with \( \mathbf{Er}(\varphi_{A,b}) \). Consequently, the following result on the error bound modulus of \( \varphi_{A,b} \) follows directly from \cref{lem3.1}.

\medskip

\begin{proposition}\label{pro3.1}
Let \( A \) be an \( m \times n \) matrix and \( b \in \mathbb{R}^m \). Then the system \( Ax \leq b \) admits a global error bound if and only if
\begin{equation}\label{3.8}
\min_{\|h\|=1}\max_{i\in J} a_i^T h < 0
\end{equation}
holds for every \( J \in \mathcal{J}_{A,b} \).
\end{proposition}

\begin{proof}
By \cref{lem3.1}, the system \( Ax \leq b \) admits a global error bound if and only if
\[
\inf_{\varphi_{A,b}(x)>0}
\left(
-\inf_{\|h\|=1} d^+\varphi_{A,b}(x,h)
\right)
>0.
\]

For any \( x \in \mathbb{R}^n \) with \( \varphi_{A,b}(x)>0 \), it follows from \cite[Theorem 4.2.3]{SIAM13} that
\[
d^+\varphi_{A,b}(x,h)
=
\max_{i\in J_{A,b}(x)} a_i^T h.
\]
Hence,
\[
-\inf_{\|h\|=1} d^+\varphi_{A,b}(x,h)
=
-\min_{\|h\|=1}\max_{i\in J_{A,b}(x)} a_i^T h.
\]
Therefore, the global error bound condition is equivalent to requiring that
\[
\min_{\|h\|=1}\max_{i\in J} a_i^T h < 0
\quad \text{for all } J \in \mathcal{J}_{A,b},
\]
which yields the desired result.\qed
\end{proof}

\vskip2mm

We are now in a position to study the computational complexity of error bounds for systems of linear inequalities. To this end, we formalize the problem as follows:

\begin{itemize}
\item[] \textbf{ERROR BOUND FOR LINEAR INEQUALITIES}
\item[] \textbf{INSTANCE:} Given positive integers \( m,n \), an \( m\times n \) matrix \( A \), and a vector \( b \in \mathbb{R}^m \), decide whether the system \( Ax \leq b \) admits a global error bound.
\end{itemize}

According to \cref{pro3.1}, determining whether the system \( Ax \leq b \) admits a global error bound is equivalent to verifying a finite family of min--max conditions of the form \eqref{3.8} and checking that all of them are strictly negative. Consequently, the complexity of error bounds is naturally reduced to the complexity of such min--max optimization problems. By \cref{lem3.2}, each individual min--max problem of this type can be solved in polynomial time.

Our main result shows that, despite this tractability at the level of individual subproblems, the overall problem is intractable in the worst case.

\begin{theorem}\label{th3.1}
\textbf{ERROR BOUND FOR LINEAR INEQUALITIES} is not in the class {\bf P}.
\end{theorem}

\begin{proof}
Suppose, by contradiction, that there exists a polynomial-time algorithm \( \Phi \) for the ERROR BOUND FOR LINEAR INEQUALITIES problem (with respect to the input size). Then, by \cref{pro3.1}, \( \Phi \) can determine whether condition \eqref{3.8} holds for all \( J \in \mathcal{J}_{A,b} \), given any instance \( (A,b) \).

We construct a family of instances for which the cardinality of \( \mathcal{J}_{A,b} \) is exponential in \( m \), which contradicts the existence of such a polynomial-time algorithm.

Assume that \( m \leq n \) and that the vectors \( \{a_1^T,\ldots,a_m^T\} \) are linearly independent. We claim that
\begin{equation}\label{3.10}
\mathcal{J}_{A,b}
=
\{J \subseteq \{1,\ldots,m\} : J \neq \emptyset\}.
\end{equation}
Hence, \( |\mathcal{J}_{A,b}| = 2^m - 1 \), which is exponential in \( m \).

It remains to prove the claim. It suffices to show that every nonempty subset \( J \subseteq \{1,\ldots,m\} \) belongs to \( \mathcal{J}_{A,b} \). Fix such a set \( J \) and let \( \varepsilon > 0 \). Define
\[
\tilde b_i := b_i + \varepsilon \quad \text{for } i \in J,
\qquad
\tilde b_i := b_i \quad \text{for } i \notin J.
\]
Consider the linear system
\[
\left\{
\begin{aligned}
a_i^T x &= \tilde b_i, \quad i \in J,\\
a_i^T x &= b_i, \quad i \notin J.
\end{aligned}
\right.
\]
Since \( \operatorname{rank}(A)=m \leq n \), this system admits a solution \( \bar x \). By construction,
\[
\varphi_{A,b}(\bar x)=\varepsilon,
\qquad
J_{A,b}(\bar x)=J,
\]
which implies \( J \in \mathcal{J}_{A,b} \). This proves \eqref{3.10}.
\qed
The exponential size of \( \mathcal{J}_{A,b} \) contradicts the existence of a polynomial-time algorithm \( \Phi \), completing the proof.
\end{proof}
\vskip2mm

\begin{remark}\label{rem3.1}

\cref{th3.1} shows that the \textbf{ERROR BOUND FOR LINEAR INEQUALITIES} problem does not belong to the class {\bf P}. However, this does not imply that the problem lies in {\bf NP}, since establishing such a membership would require resolving the well-known open problem \( \mathbf{P} \neq \mathbf{NP} \) in complexity theory.

To the best of our knowledge, there is a widespread belief that \( \mathbf{P} \neq \mathbf{NP} \), although no proof is currently available. In the absence of a resolution, it is customary to proceed under this working hypothesis rather than to attempt to prove the contrary. Accordingly, in this paper we adopt a tentative viewpoint on the structure of {\bf NP}, under the expectation (though not the certainty) that the set \( \mathbf{NP} \setminus \mathbf{P} \) is nonempty.
\end{remark}

%
%
%

Based on \cref{rem3.1}, it becomes necessary to study the complexity of the   ERROR BOUND FOR LINEAR INEQUALITIES problem within the framework of the {\bf NP} world.
\medskip

A decision problem related to   ERROR BOUND FOR LINEAR INEQUALITIES is described as follows:

\begin{itemize}
    \item[] { ERROR BOUND FOR LINEAR INEQUALITIES}
    \item[] INSTANCE: {\it Two positive integers \( m, n \), an \( m \times n \) matrix \( A \), and an \( m \)-vector \( b \).}
    \item[] QUESTION: {\it Does the system \( Ax \leq b \) admit an error bound?}
\end{itemize}

The following theorem establishes that the decision version of the   ERROR BOUND FOR LINEAR INEQUALITIES problem is Co-{\bf NP}-complete.

\begin{theorem}\label{th3.2}
The decision problem \textbf{ERROR BOUND FOR LINEAR INEQUALITIES} is Co-\textbf{NP}-complete.
\end{theorem}

\begin{proof}
We first show that if an instance has a ``no'' answer, then there exists a certificate can be verified in polynomial time.

Let \( A \) be an \( m \times n \) matrix and \( b \in \mathbb{R}^m \). Take \( \bar{x} \in \mathbb{R}^n \). Then the value \( \varphi_{A,b}(\bar{x}) \) and the active index set \( J_{A,b}(\bar{x}) \) can be computed in polynomial time. Moreover, by \cref{coro3.1}, one can decide in polynomial time whether
\[
\min_{\|h\|=1}\max_{i\in J_{A,b}(\bar{x})} a_i^T h \ge 0.
\]
Combining this with \cref{pro3.1}, it follows that if the system \( Ax \le b \) does not admit a global error bound, then \( \bar{x} \) serves as a certificate verifiable in polynomial time. Hence, the complement of the problem belongs to NP.

It remains to prove Co-\textbf{NP}-completeness by reducing SAT to the complement of such problem.

\medskip

Let \( U=\{u_1,\ldots,u_p\} \) be a set of Boolean variables and \( \mathcal{C}=\{C_1,\ldots,C_q\} \) a SAT instance. We construct a matrix \( A \) and vector \( b \) such that \( \mathcal{C} \) is satisfiable if and only if the system \( Ax \le b \) fails to admit an error bound.

We now describe the construction. Let \( n>p \) and select a maximal linearly independent subset of \( U \), denoted (without loss of generality) by \( \{u_1,\ldots,u_m\} \). Choose an arbitrary vector \( (b_1,\ldots,b_m)^T \in \mathbb{R}^m \), and let \( U^1 \subseteq U \) be a subset of minimal cardinality such that
\begin{equation}\label{3.012}
U^1 \cap C_k \neq \emptyset, \qquad k=1,\ldots,q.
\end{equation}
Define
\[
\mathcal{M}_1 := \{ i : u_i \in U^1 \}, \qquad
\mathcal{M}_2 := \{1,\ldots,m\} \setminus \mathcal{M}_1.
\]
For each \( i \in \mathcal{M}_2 \), set \( a_i := u_i \), and consider the system
\[
a_i^T x \le b_i, \qquad i \in \mathcal{M}_2,
\]
which defines \( Ax \le b \).

Given a truth assignment \( \Psi:U\to\{T,F\} \), define:
\begin{itemize}
\item[(i)] For \( j \notin \mathcal{M}_2 \), set \( \Psi(u_j)=T \) iff
\[
\min_{\|h\|=1}\max_{i\in\mathcal{M}_2} a_i^T h \ge 0.
\]

\item[(ii)] For \( j \in \mathcal{M}_2 \), set \( \Psi(u_j)=T \) iff
\[
\min_{\|h\|=1}\max_{i\in\mathcal{M}_2\setminus\{j\}} a_i^T h \ge 0.
\]
\end{itemize}

We now prove that this construction yields a valid reduction, i.e., \( \mathcal{C} \) is satisfiable if and only if \( Ax \le b \) does not admit an error bound.

\medskip

\noindent\textbf{Necessity.}
Assume that \( \mathcal{C} \) is satisfiable. Then for each clause \( C_k \), there exists \( u_j \in C_k \) with \( \Psi(u_j)=T \), and hence
\begin{equation}\label{3.12}
\min_{\|h\|=1}\max_{i\in\mathcal{M}_2} a_i^T h \ge 0.
\end{equation}

Let \( \varepsilon>0 \), and consider the perturbed system
\begin{equation}\label{3.13}
a_i^T x = b_i + \varepsilon, \qquad i \in \mathcal{M}_2.
\end{equation}
Since the vectors \( \{a_i : i\in\mathcal{M}_2\} \) are linearly independent, system \eqref{3.13} admits a solution \( \hat{x} \). For this point,
\[
\varphi_{A,b}(\hat{x})=\varepsilon>0,
\qquad
J_{A,b}(\hat{x})=\mathcal{M}_2.
\]
By \cref{pro3.1}, it follows that \( Ax \le b \) does not admit an error bound.

\medskip

\noindent\textbf{Sufficiency.}
Assume that \( Ax \le b \) does not admit an error bound. By \cref{pro3.1}, there exists \( \tilde{x} \in \mathbb{R}^n \) such that
\[
\varphi_{A,b}(\tilde{x})>0,
\qquad
\min_{\|h\|=1}\max_{i\in J_{A,b}(\tilde{x})} a_i^T h \ge 0.
\]
Since \( J_{A,b}(\tilde{x}) \subseteq \mathcal{M}_2 \), we obtain
\[
\min_{\|h\|=1}\max_{i\in\mathcal{M}_2} a_i^T h \ge 0.
\]

By \eqref{3.012}, each clause \( C_k \) contains some \( u_j \in U^1 \cap C_k \). Since \( u_j \notin \mathcal{M}_2 \), the definition of \( \Psi \) implies \( \Psi(u_j)=T \). Hence every clause is satisfied, and \( \mathcal{C} \) is satisfiable.\qed
\end{proof}

\medskip

We are inspired by \cref{th3.2} to study the complement of the   ERROR BOUND FOR LINEAR INEQUALITIES problem and its complexity. A decision problem related to the complement is described as follows:

\begin{itemize}
    \item[] {\bf Complement of   ERROR BOUND FOR LINEAR INEQUALITIES}
    \item[] INSTANCE: {\it Two positive integers \( m, n \), an \( m \times n \) matrix \( A \), and an \( m \)-vector \( b \).}
    \item[] QUESTION: {\it Does the system \( Ax \leq b \) admit no error bound?}
\end{itemize}

\medskip

The following corollary, which follows immediately from \cref{th3.2} and \cref{coro3.1}, characterizes the complexity of the complement of   ERROR BOUND FOR LINEAR INEQUALITIES.

\begin{corollary}\label{coro3.3}
The decision problem Complement of   ERROR BOUND FOR LINEAR INEQUALITIES is {\bf NP}-complete and therefore {\bf NP}-hard.
\end{corollary}

It follows from \cref{th3.1} and \cref{th3.2} that {``numbers"} 
play a significant role in the complexity proofs. We now focus our attention on analyzing the complexity of the complement of ERROR BOUND FOR LINEAR INEQUALITIES  from this numerical perspective. To this end, we recall the concepts of {\it pseudo-polynomial time algorithms} and {\it strong NP-completeness}.

\medskip

Given a problem $\Pi$, let $\mathcal{D}_{\Pi}$ denote the set of its instances. 
	Assume that $\Pi$ is equipped with two encoding-independent functions
	\[
	\textit{Length} : \mathcal{D}_{\Pi} \to \mathbb{N}, 
	\qquad 
	\textit{Max} : \mathcal{D}_{\Pi} \to \mathbb{N}.
	\]
	Here, $\textit{Length}[\Lambda]$ denotes the number of symbols required to encode 
	an instance $\Lambda \in \mathcal{D}_{\Pi}$, and 
	$\textit{Max}[\Lambda]$ denotes the largest magnitude of the integers appearing in $\Lambda$, 
	taken to be $0$ if no integers occur.

Recall from \cite{GJ79} that an algorithm solving $\Pi$ is called 
	\emph{pseudo-polynomial} if its running time is bounded by a polynomial 
	in both $\textit{Length}[\Lambda]$ and $\textit{Max}[\Lambda]$.

A problem $\Pi$ is called a \emph{number problem} if there exists no polynomial $P(\cdot)$ 
	such that $\textit{Max}[\Lambda] \le P(\textit{Length}[\Lambda])$ for all 
	$\Lambda \in \mathcal{D}_{\Pi}$. 
	For any polynomial $P(\cdot)$, let $\Pi_P$ denote the restriction of $\Pi$ 
	to instances satisfying $\textit{Max}[\Lambda] \le P(\textit{Length}[\Lambda])$. 
	We say that $\Pi$ is \emph{strongly NP-complete} if $\Pi \in \mathbf{NP}$ 
	and there exists a polynomial $P(\cdot)$ such that $\Pi_P$ is NP-complete.

\begin{theorem}\label{th3.3}
There exists a pseudo-polynomial time algorithm to solve the complement of ERROR BOUND FOR LINEAR INEQUALITIES.
\end{theorem}

{\bf Proof.} Let \( \Lambda \) denote the instance:  ``Given an \( m \times n \) matrix \( A \) and an \( m \)-vector \( b \), does the system of linear inequalities \( Ax \leq b \) admit an error bound?" 

Denote by \( \textit{Length}[\Lambda] \) the  encoding length of the input and by \( \textit{Max}[\Lambda] \) the largest integer appearing in \( \Lambda \).

Define
\begin{equation}\label{3.18-260328}
  \widehat{\mathcal{J}}_{A,b}:=\left\{J\in \mathcal{J} _{A,b}:	\min_{\|h\| = 1} \max_{i \in J} a_i^T h \geq 0\right\}.
\end{equation}
By \cref{pro3.1}, the system $Ax\leq b$ does not admit an error bound if and only if $  \widehat{\mathcal{J}}_{A,b}\not=\emptyset$. 

Let $\hat J\in  \widehat{\mathcal{J}}_{A,b}$. Then
\[
\big |\hat J \big | \leq m \leq \textit{Max}[\Lambda].
\]
Define
\[
\mathcal{J}_1 :=\left \{ J \in   \mathcal{J}_{A,b}: |J| = |\hat J|  \right\}.
\]
Observe that the cardinality of \( \mathcal{J}_1 \) is bounded by \( \binom{m}{|\hat{J}|} \), and hence is polynomially bounded in \( \textit{Max}[\Lambda] \) when \( |\hat{J}| \) is fixed.

We claim that the system \( Ax \leq b \) does not admit an error bound if and only if there exists \( J_1 \in \mathcal{J}_1 \) such that
\begin{equation}\label{3.15}
	\min_{\|h\| = 1} \max_{i \in J_1} a_i^T h \geq 0.
\end{equation}
Assuming the claim, it follows from \cref{coro3.1} that the above condition can be verified in polynomial time, which yields the desired algorithm for the complement problem.

It remains to prove the claim.

\medskip
\noindent\textbf{Necessity.}
If \( Ax \leq b \) does not admit an error bound, then there exists 
\( \hat{J} \in \widehat{\mathcal{J}}_{A,b} \). By construction, \( \hat{J} \in \mathcal{J}_1 \), and \eqref{3.15} holds.

\medskip
\noindent\textbf{Sufficiency.}
Conversely, suppose that there exists \( J_1 \in \mathcal{J}_1 \) satisfying \eqref{3.15}. Since \( \mathcal{J}_1 \subseteq \mathcal{J}_{A,b} \), it follows from \cref{pro3.1} that the system \( Ax \leq b \) does not admit an error bound.

Hence the claim holds. \qed
\medskip

Based on \cref{th3.3}, the following observation follows from \cite[Observation 4.1 and Observation 4.2]{GJ79}, under the hypothesis that \( \mathbf{P} \neq \mathbf{NP} \).

\medskip

\noindent{\bf Observation 3.1.} {\it The complement of  ERROR BOUND FOR LINEAR INEQUALITIES is a number problem and is not {\bf NP}-complete in the strong sense unless \( \mathbf{P} = \mathbf{NP} \).}

\medskip

\section{Concluding Remarks and Perspectives}
It has been more than seventy years since Hoffman's pioneering work on error bounds for systems of linear inequalities. In contrast to most existing studies, which primarily focus on establishing sufficient or necessary conditions for the validity of error bounds, this paper investigates the computational complexity of error bounds in linear inequality systems.

From the perspective of practical applications, understanding the computational complexity of error bounds opens a new avenue for research. For example, knowing that the recognition problem is Co-NP-complete provides valuable insight into which research directions are likely to be productive. One possible direction is the identification of efficient algorithms for special classes of linear systems. Another is the design of algorithms that, although not polynomial-time in the worst case, perform efficiently on instances arising in practice. Alternatively, one may consider relaxed or approximate formulations that are more amenable to efficient computation.

In this sense, the complexity results established in this paper are not only of theoretical interest but also provide guidance for future research. They help clarify which approaches are most likely to lead to computationally efficient methods for detecting error bounds and analyzing their stability under data perturbations.

Ultimately, the study of the complexity of error bound recognition is intended to direct research toward methods with the greatest potential for computational feasibility and practical effectiveness.

At the same time, the techniques developed here for analyzing the complexity of error bounds may also be applicable to other fundamental notions and properties in variational analysis and optimization. A promising direction for future work is to undertake preliminary complexity investigations for special problems related to error bounds--some potentially more difficult, such as stability-related questions, and others possibly more tractable--before pursuing a deeper and more comprehensive complexity analysis. This line of inquiry forms part of our future research agenda.

Several open questions remain concerning the complexity of error bounds for systems of linear inequalities:

\begin{itemize}
    \item[(i)] Is the ERROR BOUND FOR LINEAR INEQUALITIES problem {\bf NP}- complete?
    \item[(ii)] Does there exist a pseudo-polynomial time algorithm to solve the   
    ERROR BOUND FOR LINEAR INEQUALITIES  problem?
    \item[(iii)] What is the complexity of stability of error bounds for systems of linear inequalities?
\end{itemize}

\section*{ Acknowledgments}
We would like to thank Professor Claudia  Sagastiz\'abal  for her valuable comments and suggestions, which have helped improve the presentation of this work. We are also grateful to Dr. Nguyen Duy Cuong for kindly providing his PhD dissertation, ``Transversality, Regularity and Error Bounds in Variational Analysis". Finally, we thank the two anonymous reviewers for their careful reading and detailed comments, which have contributed to improving the paper.

\medskip

\section*{Data Availability} Data sharing not applicable to this article as no datasets were generated or analysed during
the current study.

\section* {Conflicts of Interest: }

The authors declare that they have no conflicts of interest, and the manuscript has no associated data.


\begin{thebibliography}{10}
	
	\bibitem{Abassi-Thera1}
	M.~Abbasi and M.~Th\'{e}ra.
	\newblock Strongly regular points of mappings.
	\newblock {\em Fixed Point Theory Algorithms Sci. Eng.}, pages Paper No. 14, 13
	pp., 2021.
	
	\bibitem{Abassi-Thera}
	M.~Abbasi and M.~Th\'{e}ra.
	\newblock About error bounds in metrizable topological vector spaces.
	\newblock {\em Set-Valued Var. Anal.}, 30(4):1291--1311, 2022.
	
	\bibitem{AC1988}
	A.~Auslender and J.-P. Crouzeix.
	\newblock Global regularity theorems.
	\newblock {\em Math. Oper. Res.}, 13(2):243--253, 1988.
	
	\bibitem{AF96}
	D.~Avis and K.~Fukuda.
	\newblock Reverse search for enumeration.
	\newblock {\em Discrete Appl. Math.}, 65:21--46, 1996.
	
	\bibitem{5}
	H.~H. Bauschke and J.~M. Borwein.
	\newblock On projection algorithms for solving convex feasibility problems.
	\newblock {\em SIAM Rev.}, 38(3):367--426, 1996.
	
	\bibitem{6}
	A.~Beck and M.~Teboulle.
	\newblock Convergence rate analysis and error bounds for projection algorithms
	in convex feasibility problems.
	\newblock {\em Optim. Methods Softw.}, 18(4):377--394, 2003.
	
	\bibitem{7}
	E.~M. Bednarczuk and A.~Y. Kruger.
	\newblock Error bounds for vector-valued functions: necessary and sufficient
	conditions.
	\newblock {\em Nonlinear Anal.}, 75(3):1124--1140, 2012.
	
	\bibitem{BurDeng02}
	J.~V. Burke and S.~Deng.
	\newblock Weak sharp minima revisited. i. basic theory.
	\newblock {\em Control Cybernet.}, 31(3):439--469, 2002.
	
	\bibitem{BD2}
	J.~V. Burke and S.~Deng.
	\newblock Weak sharp minima revisited. ii. application to linear regularity and
	error bounds.
	\newblock {\em Math. Program.}, 104(2--3):235--261, 2005.
	
	\bibitem{CKLT}
	M.~J. C\'{a}novas, A.~Y. Kruger, M.~A. L\'{o}pez, J.~Parra, and M.~A.
	Th\'{e}ra.
	\newblock Calmness modulus of linear semi-infinite programs.
	\newblock {\em SIAM J. Optim.}, 24(1):29--48, 2014.
	
	\bibitem{CK70}
	D.~R. Chand and S.~S. Kapur.
	\newblock An algorithm for convex polytopes.
	\newblock {\em J. ACM}, 17:78--86, 1970.
	
	\bibitem{CCH53}
	A.~Charnes, W.~W. Cooper, and A.~Henderson.
	\newblock {\em An Introduction to Linear Programming}.
	\newblock Wiley, 1953.
	
	\bibitem{Chazelle1993}
	B.~Chazelle.
	\newblock An optimal convex hull algorithm in any fixed dimension.
	\newblock {\em Discrete Comput. Geom.}, 10:377--409, 1993.
	
	\bibitem{16}
	P.~L. Combettes.
	\newblock Hilbertian convex feasibility problem: convergence of projection
	methods.
	\newblock {\em Appl. Math. Optim.}, 35(3):311--330, 1997.
	
	\bibitem{CCZ2014}
	M.~Conforti, G.~Cornu\'ejols, and G.~Zambelli.
	\newblock {\em Integer Programming}.
	\newblock Springer, Switzerland, 2014.
	
	\bibitem{Cook71}
	S.~Cook.
	\newblock The complexity of theorem proving procedures.
	\newblock In {\em Proceedings of the Third Annual ACM Symposium on Theory of
		Computing}, pages 151--158, 1971.
	
	\bibitem{CK2022}
	N.~D. Cuong and A.~Y. Kruger.
	\newblock Error bounds revisited.
	\newblock {\em Optimization}, 71(4):1021--1053, 2022.
	
	\bibitem{DuMa2021}
	J.~Dutta and J.~E. Mart\'{i}nez-Legaz.
	\newblock Error bounds for inequality systems defining convex sets.
	\newblock {\em Math. Program.}, 189(1--2):1--25, 2021.
	
	\bibitem{Dye83}
	M.~E. Dyer.
	\newblock The complexity of vertex enumeration methods.
	\newblock {\em Math. Oper. Res.}, 8:381--402, 1983.
	
	\bibitem{FHKO2010}
	M.~Fabian, R.~Henrion, A.~Y. Kruger, and J.~Outrata.
	\newblock Error bounds: necessary and sufficient conditions.
	\newblock {\em Set-Valued Var. Anal.}, 18(2):121--149, 2010.
	
	\bibitem{Fukuda2004}
	K.~Fukuda.
	\newblock Frequently asked questions in polyhedral computation, 2004.
	\newblock ETH Z\"{u}rich.
	
	\bibitem{GJ79}
	M.~Garey and D.~Johnson.
	\newblock {\em Computers and Intractability: A Guide to the Theory of
		NP-Completeness}.
	\newblock Freeman, 1979.
	
	\bibitem{Hoffman1952}
	A.~J. Hoffman.
	\newblock On approximate solutions of systems of linear inequalities.
	\newblock {\em J. Res. Natl. Bur. Stand.}, 49:263--265, 1952.
	
	\bibitem{Ioffe1979}
	A.~D. Ioffe.
	\newblock Regular points of lipschitz functions.
	\newblock {\em Trans. Amer. Math. Soc.}, 251:61--69, 1979.
	
	\bibitem{ioffe-JAMS-1}
	A.~D. Ioffe.
	\newblock Metric regularity: a survey. part i. theory.
	\newblock {\em J. Aust. Math. Soc.}, 101(2):188--243, 2016.
	
	\bibitem{ioffe-JAMS-2}
	A.~D. Ioffe.
	\newblock Metric regularity: a survey. part ii. applications.
	\newblock {\em J. Aust. Math. Soc.}, 101(3):376--417, 2016.
	
	\bibitem{ioffe-book}
	A.~D. Ioffe.
	\newblock {\em Variational Analysis of Regular Mappings}.
	\newblock Springer, Cham, 2017.
	
	\bibitem{SIAM13}
	A.~D. Ioffe and V.~M. Tikhomirov.
	\newblock {\em Theory of Extremal Problems}.
	\newblock North-Holland, Amsterdam, 1979.
	
	\bibitem{jourani}
	A.~Jourani.
	\newblock Hoffman's error bound, local controllability, and sensitivity
	analysis.
	\newblock {\em SIAM J. Control Optim.}, 38(3):947--970, 2000.
	
	\bibitem{Khachiyan1995}
	L.~Khachiyan.
	\newblock On the complexity of approximating extremal determinants in matrices.
	\newblock {\em J. Complexity}, 11(1):138--153, 1995.
	
	\bibitem{KlatteLi1999}
	D.~Klatte and W.~Li.
	\newblock Asymptotic constraint qualifications and global error bounds for
	convex inequalities.
	\newblock {\em Math. Program.}, 84(1):137--160, 1999.
	
	\bibitem{Kru15.2}
	A.~Y. Kruger.
	\newblock Error bounds and h\"older metric subregularity.
	\newblock {\em Set-Valued Var. Anal.}, 23(4):705--736, 2015.
	
	\bibitem{Kru15}
	A.~Y. Kruger.
	\newblock Error bounds and metric subregularity.
	\newblock {\em Optimization}, 64(1):49--79, 2015.
	
	\bibitem{Kruger-LY}
	A.~Y. Kruger, M.~A. L\'opez, X.~Yang, and J.~Zhu.
	\newblock H\"older error bounds and h\"older calmness with applications to
	convex semi-infinite optimization.
	\newblock {\em Set-Valued Var. Anal.}, 27(4):995--1023, 2019.
	
	\bibitem{KLT2017}
	A.~Y. Kruger, D.~R. Luke, and N.~H. Thao.
	\newblock Subtransversality of collections of sets.
	\newblock {\em Set-Valued Var. Anal.}, 25(4):701--729, 2017.
	
	\bibitem{LewisPang1996}
	A.~S. Lewis and J.-S. Pang.
	\newblock Error bounds for convex inequality systems.
	\newblock In {\em Generalized Convexity, Generalized Monotonicity}, volume~27
	of {\em Nonconvex Optim. Appl.}, pages 75--110. Kluwer Acad. Publ.,
	Dordrecht, 1998.
	
	\bibitem{LMY2018}
	M.~H. Li, K.~W. Meng, and X.~Q. Yang.
	\newblock On error bound moduli for locally lipschitz and regular functions.
	\newblock {\em Math. Program.}, 171(1):463--487, 2018.
	
	\bibitem{HLu}
	D.~R. Luke.
	\newblock Nonconvex notions of regularity and convergence of fundamental
	algorithms for feasibility problems.
	\newblock {\em SIAM J. Optim.}, 23(4):2397--2419, 2013.
	
	\bibitem{Mangasarian1985}
	O.~L. Mangasarian.
	\newblock A condition number for differentiable convex inequalities.
	\newblock {\em Math. Oper. Res.}, 10(2):175--179, 1985.
	
	\bibitem{Boris1}
	B.~S. Mordukhovich.
	\newblock {\em Variational Analysis and Generalized Differentiation I}.
	\newblock Springer, Berlin, 2006.
	
	\bibitem{Boris2}
	B.~S. Mordukhovich.
	\newblock {\em Variational Analysis and Generalized Differentiation II}.
	\newblock Springer, Berlin, 2006.
	
	\bibitem{NT2004}
	H.~V. Ngai and M.~Th\'era.
	\newblock Error bounds and implicit multifunction theorem in smooth banach
	spaces and applications to optimization.
	\newblock {\em Set-Valued Anal.}, 12(1--2):195--223, 2004.
	
	\bibitem{Huynh2008Error}
	H.~V. Ngai and M.~Th\'era.
	\newblock Error bounds in metric spaces and application to the perturbation
	stability of metric regularity.
	\newblock {\em SIAM J. Optim.}, 19(1):1--20, 2008.
	
	\bibitem{NT2009}
	H.~V. Ngai and M.~Th\'era.
	\newblock Error bounds for systems of lower semicontinuous functions in asplund
	spaces.
	\newblock {\em Math. Program.}, 116(1--2):397--427, 2009.
	
	\bibitem{Pang1997}
	J.-S. Pang.
	\newblock Error bounds in mathematical programming.
	\newblock {\em Math. Program.}, 79(1--3):299--332, 1997.
	
	\bibitem{penot-book}
	J.-P. Penot.
	\newblock {\em Calculus without Derivatives}.
	\newblock Springer, New York, 2013.
	
	\bibitem{Ph}
	R.~R. Phelps.
	\newblock {\em Convex Functions, Monotone Operators and Differentiability}.
	\newblock Springer, Berlin, 2 edition, 1993.
	
	\bibitem{Rob73}
	S.~M. Robinson.
	\newblock Bounds for error in the solution set of a perturbed linear program.
	\newblock {\em Linear Algebra Appl.}, 6:69--81, 1973.
	
	\bibitem{Robinson1975}
	S.~M. Robinson.
	\newblock An application of error bounds for convex programming in a linear
	space.
	\newblock {\em SIAM J. Control}, 13:271--273, 1975.
	
	\bibitem{Rob77}
	S.~M. Robinson.
	\newblock A characterization of stability in linear programming.
	\newblock {\em Oper. Res.}, 25(3):435--447, 1977.
	
	\bibitem{thibault-book}
	L.~Thibault.
	\newblock {\em Unilateral Variational Analysis in Banach Spaces. Part I:
		General Theory}.
	\newblock World Scientific, Singapore, 2023.
	
	\bibitem{TsB93}
	P.~Tseng and D.~P. Bertsekas.
	\newblock On the convergence of the exponential multiplier method for convex
	programming.
	\newblock {\em Math. Program.}, 60(1):1--19, 1993.
	
	\bibitem{Tuncel1999}
	L.~Tuncel.
	\newblock Approximating the complexity measure of the vavasis--ye algorithm is
	{NP}-hard.
	\newblock {\em Math. Program.}, 86(1):219--223, 1999.
	
	\bibitem{WTY2022}
	Z.~Wei, M.~Th\'era, and J.-C. Yao.
	\newblock Characterizations of stability of error bounds for convex inequality
	constraint systems.
	\newblock {\em Open J. Math. Optim.}, 3:Art. No. 2, 17 pp., 2022.
	
	\bibitem{WTY2023}
	Z.~Wei, M.~Th\'era, and J.-C. Yao.
	\newblock Primal characterizations of error bounds for composite-convex
	inequalities.
	\newblock {\em J. Convex Anal.}, 30(4):1329--1350, 2023.
	
	\bibitem{WTY2024}
	Z.~Wei, M.~Th\'era, and J.-C. Yao.
	\newblock Primal characterizations of stability of error bounds for
	semi-infinite convex constraint systems.
	\newblock {\em Optimization}, 73(12):3725--3753, 2024.
	
	\bibitem{WTY2024-SVVA}
	Z.~Wei, M.~Th\'era, and J.-C. Yao.
	\newblock Subtransversality and strong {CHIP} of closed sets in asplund spaces.
	\newblock {\em Set-Valued Var. Anal.}, 32:23, 2024.
	
	\bibitem{WTY2025}
	Z.~Wei, M.~Th\'era, and J.-C. Yao.
	\newblock Perturbation analysis of error bounds for convex functions on banach
	spaces.
	\newblock {\em J. Convex Anal.}, 32(3):883--900, 2025.
	
	\bibitem{WZ}
	Z.~Wei, J.-C. Yao, and X.~Y. Zheng.
	\newblock Strong abadie {CQ}, {ACQ}, calmness and linear regularity.
	\newblock {\em Math. Program.}, 145:97--131, 2014.
	
	\bibitem{55}
	Z.~Wu and J.~J. Ye.
	\newblock On error bounds for lower semicontinuous functions.
	\newblock {\em Math. Program.}, 92(2):301--314, 2002.
	
	\bibitem{ZN2004}
	X.~Y. Zheng and K.~F. Ng.
	\newblock Metric regularity and constraint qualifications for convex
	inequalities on banach spaces.
	\newblock {\em SIAM J. Optim.}, 14(3):757--772, 2004.
	
	\bibitem{Za2001}
	C.~Z\u{a}linescu.
	\newblock Weak sharp minima, well-behaving functions and global error bounds
	for convex inequalities in banach spaces.
	\newblock In {\em Proceedings of the 12th Baikal International Conference on
		Optimization Methods and Their Applications}, pages 272--284, Irkutsk,
	Russia, 2001.
	
\end{thebibliography}
\end{document}